\documentclass{birkjour}
\usepackage[OT1]{fontenc}
\usepackage{amsthm,amsmath,amsfonts,amssymb}
\usepackage{amsmath,amsfonts,amssymb}
\usepackage[colorlinks]{hyperref}

\theoremstyle{plain}% default
\newtheorem{thm}{Theorem}[section]

\newtheorem{prop}[thm]{Proposition}
\newtheorem{cor}[thm]{Corollary}

\theoremstyle{definition}

\newtheorem{conj}[thm]{Conjecture}
\newtheorem{exmp}[thm]{Example}

\theoremstyle{remark}
\newtheorem{rem}[thm]{Remark}

\newcommand{\indep}{{\;\bot\!\!\!\!\!\!\bot\;}}

\begin{document}

\title{Binary cumulant varieties}
\author{Bernd Sturmfels}
\address{Department of Mathematics, University
of California, Berkeley, CA 94720, USA, {\tt bernd@math.berkeley.edu}}

\author{Piotr Zwiernik}
\address{Institute Mittag-Leffler, 
18260 Djursholm, Sweden,
{\tt piotr.zwiernik@gmail.com}}
\thanks{This work was conducted during the Spring 2011 program
{\em ``Algebraic Geometry with a View Towards Applications''} at
the Institute Mittag-Leffler in Djursholm, Sweden. The first author was 
partially supported by KTH Stockholm and the US National Science Foundation (DMS-0968882). The second author gratefully acknowledges support from the AXA Mittag-Leffler Fellowship Project, sponsored by the AXA Research Fund.}

\markboth{authors}{Algebraic geometry of binary cumulants}
\
 \subjclass{Primary: 13P25; Secondary: 05A40, 14Q15, 60C05}
\keywords{Algebraic statistics, cumulants, moments, binary data, context-specific independence,
hyperdeterminant, Segre variety, secant variety}
\date{}

\begin{abstract}

Algebraic statistics for binary random variables is concerned with
highly structured algebraic varieties in the space of 
$2 {\times} 2 {\times} \cdots {\times} 2$-tensors. We demonstrate the advantages
of representing such varieties in the coordinate system of binary cumulants.
Our primary focus lies on hidden subset models.
Parametrizations and implicit equations in cumulants are derived for
hyperdeterminants, for secant and tangential varieties of Segre varieties,
and for certain context-specific independence models.
Extending work of Rota and collaborators, we explore the 
polynomial inequalities satisfied by cumulants.

\end{abstract}
%Version: \today
\maketitle
\section{Introduction}\label{sec:introduction}

Cumulants have a long and interesting history dating back to Thorvald N. Thiele, a Danish mathematician, who introduced them in 1889.
See \cite{hald2000early} for a historical perspective.
  The main motivation to study them was that multivariate probability distributions are often easier to analyze when expressed in terms of cumulants. Moreover, cumulants are especially useful when dealing with the normal distribution, and hence they are a critical tool in asymptotic statistics (see e.g. \cite{BrWy, mccullagh1987tms,
speed1983cumulants, torney}). 
 Various invariance properties of cumulants make them  
interesting  also from an algebraic or combinatorial point of view.
 Rota and his collaborators
\cite{bruno1999probability, rota_cumulants} developed a
combinatorial theory of cumulants, and, more recently, 
Pistone and Wynn introduced {\em cumulant varieties}  \cite{pistonewynn}  into algebraic statistics. These concepts gave rise to umbral calculus  \cite{rota1978umbral},
an approach to combinatorial sequences using cumulants.

Building on this circle of ideas, we show how cumulants can be used to study algebraic varieties
in tensor spaces. Thus, cumulants can be also used outside of the probabilistic context where we deal with sequences of nonnegative numbers summing to $1$. Here we focus on  binary states.
 Let $P=[p_{I}]_{I\subseteq [n]}$ be an $n$-dimensional  
$2{\times} 2{\times}\cdots{\times} 2$ table of complex numbers such that $\sum p_{I}=1$. We call such tensors \textit{distributions}. 
In statistical contexts one assumes in addition that the $p_I$ are real and nonnegative in which case we call them \textit{probability distributions}. In algebraic statistics, the probabilities $p_{I}$ form the
coordinates of the ambient space containing statistical models. 
For an introduction to this geometric point of view see \cite{OWbook}.

We represent the distribution $P$ by the {\em probability generating function}
$$
P(x)\quad=\quad\sum_{I\subseteq [n]} p_{I}\prod_{i\in I}x_{i}.
$$
Here $[n] = \{1,2,\ldots,n\}$ and we identify  our tables with
functions on subsets of $[n]$. In the probabilistic context we occasionally refer to the random vector $X=(X_{1},\ldots,X_{n})$ with values in $\{0,1\}^{n}$ and distribution $P$. We use here the natural identification of a subset $I\subseteq [n]$ with its support vector. 
An alternative representation of $P$ is the table of moments $M=[\mu_{I}]_{I\subseteq [n]}$, where
\begin{equation}\label{eq:mominprob}
\mu_{I} \,\,= \,\,\sum_{J\supseteq I}p_{J}.
\end{equation}
The {\em moment generating function} is a square-free polynomial in $n$ unknowns:
\begin{equation}\label{eq:mgf}
 M(x) \quad =\quad P(x_{1}+1,\ldots,x_{n}+1)\quad= \quad \sum_{I \subseteq [n]} \mu_I \prod_{i \in I} x_i . \end{equation}
The logarithm of the moment generating function gives the cumulants:
\begin{equation}
\label{eq:takelog}
 K(x) \,\,\, = \,\, \sum_{I \subseteq [n]} k_I \prod_{i \in I} x_I \quad
:= \quad {\rm log}(M(x)). 
\end{equation}
Note that $\mu_\emptyset = 1$ and $k_\emptyset = 0$. 
For the logarithm we use the familiar series 
$ {\rm log}(1 + t) 
= \sum_{i=1}^\infty (-1)^{i-1} t^i/i $. That expansion is understood
modulo the ideal $\langle x_1^2 , x_2^2 ,\ldots, x_n^2 \rangle$.
The moments can then be recovered from the cumulants via
\begin{equation}
\label{eq:takeexp}
 M(x) \quad = \quad {\rm exp}(K(x)). 
 \end{equation}

The transformations (\ref{eq:takelog}) and (\ref{eq:takeexp}) between
moments $\mu_I$ and cumulants $k_I$ can be written as explicit combinatorial formulas
(see e.g. \cite[\S 2.3]{mccullagh1987tms}, \cite{rota_cumulants, speed1983cumulants}). 
  Given any $I\subseteq [n]$, let $\Pi(I)$ be the lattice of all set partitions of $I$.  We have
  \begin{equation}\label{eq:cum1}
k_{I}\quad\, =\,\, \sum_{\pi\in\Pi(I)}(-1)^{|\pi|-1}(|\pi|-1)!\prod_{B\in \pi} \mu_{B}.
\end{equation}
The sum is over partitions of $I$, the product is over blocks of a partition, and $|\pi|$ denotes the number of blocks of $\pi$. 
The moments in terms of cumulants~are
\begin{equation}\label{eq:momincum}
\mu_{I} \quad = \,\,\sum_{\pi\in\Pi(I)}\prod_{B\in \pi} k_{B}\qquad\mbox{ for all } I\subseteq [n].
\end{equation}
For instance, $I=\{1,2,3\}$ has five partitions $123$, $1|23$, $2|13$, $12|3$, $1|2|3$, and
\begin{equation}\label{eq:kappa3}
\begin{matrix}
{k}_{123}& =& \mu_{123}-\mu_{1}\mu_{23}-\mu_{2}\mu_{13}-\mu_{12}\mu_{3}+2\mu_{1}\mu_{2}\mu_{3},\\
\mu_{123} &=& k_{123} + k_{12} k_3 + k_{13} k_2 + k_{23} k_1 + k_1 k_2 k_3. 
\end{matrix}
\end{equation}
The transformation from (\ref{eq:momincum}) to (\ref{eq:cum1})
is the {\em M\"obius inversion} on the {\em partition lattice} $\Pi([n]) $,
as seen in enumerative combinatorics \cite[Exercise 3.44]{stanley2006enumerative}. 

This article is organized as follows.
In Section 2 we study the expression of {\em hyperdeterminants} in terms of cumulants.
In Section 3 we show that $SL(2)^n$-invariant tensor varieties 
are defined by $\mathbb{Z}^n$-homogeneous polynomials in the higher order cumulants
$k_I$ with $|I | \geq 2$. Section 4 concerns secants and tangents of the Segre variety,
and we show (in Theorem \ref{thm:tangtoric})
that the {\em tangential variety} becomes toric in cumulant coordinates.
A conceptual explanation for this arises from our theory of
{\em hidden subset models}, developed in Section 5. Here the main result is
Theorem \ref{thm:bigCSImodel}.
Section 6 offers an algebraic study of the {\em context-specific independence
models} due to Georgi and Schliep \cite{Georgi15072006}.
 Section 7 explores the {\em semialgebraic constraints} on
cumulants arising from probabilities,
and it addresses a conjecture proposed in
\cite{bruno1999probability}.

\section{Hyperdeterminants}\label{sec:hyperdet}

One of the most intriguing polynomial functions on $2 {\times} 2 {\times} \cdots {\times} 2$-tables
is the {\em hyperdeterminant} ${\rm Det}(P)$, which is a generalization of the determinant of a $2\times 2$ matrix.
The hyperdeterminant, first introduced by Cayley in 1843, has many equivalent definitions (see \cite{gelfand1994dra}). One of them states that  ${\rm Det}(P)$
 is the  (unique up to scaling) irreducible polynomial  in the $p_I$ that vanishes whenever the complex
hypersurface defined by the equation $P(x) = 0$ has a singular point in $\mathbb{C}^n$. Algebraically,
 the hyperdeterminant ${\rm Det}(P)$ is obtained by eliminating the
 $n$ unknowns $x_1,x_2,\ldots,x_n$ from the $n+1$ equations
$$ P(x) \,=\, \frac{\partial P}{\partial x_1} (x) =
\frac{\partial P}{\partial x_2} (x) =
 \cdots = \frac{\partial P}{\partial x_n}(x)  = 0 . $$
According to \cite[\S 14.2]{gelfand1994dra},   ${\rm Det}(P)$ is a 
 homogeneous polynomial of degree $C_n$ in the $2^n$ unknowns, where
$\sum_{n=0}^\infty C_n z^n/n! \, = \,e^{-2x}/(1-x)^2$.
So, the degrees of our hyperdeterminants are
$\,C_2 = 2, C_3 = 4,C_4 =  24, C_5 = 128 $ etc.
 
We work in the $(2^n-1)$-dimensional affine space of distributions defined by $ \sum_I p_I = 1$, or
$\mu_\emptyset = 1$. We
seek to express the hyperdeterminant on that affine space
in terms of the cumulants $k_I$. From such an expression one recovers
a formula for ${\rm Det}(P)$ in terms of the original coordinates $p_I$,
up to scaling, by using (\ref{eq:mominprob}) and~(\ref{eq:cum1}).

If $n=2$ then 
the hyperdeterminant is the  determinant of a $2 \times 2$-matrix,
$$ P \,\,= \,\, \left[\begin{array}{cc}
p_{\emptyset} & p_{2}\\
p_{1} & p_{12}
\end{array}\right] \! . $$
In statistics, this represents the {\em independence model}
for two binary random variables, and we recover the well-known fact that independence 
is equivalent to vanishing of the covariance
$$
{\rm Det}(P)\,\,= \,\,p_{12}p_{\emptyset}-p_{1}p_{2}\,\, = \,\, \mu_{12}-\mu_{1}\mu_{2} \,\, = \,\, k_{12}.
$$
The statistical meaning of  larger hyperdeterminants
will be discussed  later.  See, in particular,
the context-specific independence model in Example \ref{ex:hypercsi}.

If $n=3$ then, by \cite[Proposition 14.1.7]{gelfand1994dra}, the hyperdeterminant equals
\begin{eqnarray*}{\rm Det}(P)\,=\,
 \mu_1^2 \mu_{23}^2
{+}\mu_2^2 \mu_{13}^2
{+}\mu_3^2 \mu_{12}^2 
{+}\mu_{123}^2  +4 (\mu_1 \mu_2 \mu_3 \mu_{123}+ \mu_{12} \mu_{13} \mu_{23} )
\qquad \qquad  \\
\! -
2 (\mu_1 \mu_2 \mu_{13} \mu_{23}
{+}\mu_1 \mu_3 \mu_{12} \mu_{23} 
{+}\mu_2 \mu_3 \mu_{12} \mu_{13}  {+}
\mu_1 \mu_{23} \mu_{123} {+}
 \mu_2 \mu_{13} \mu_{123} {+} \mu_3 \mu_{12} \mu_{123} ) 
\end{eqnarray*}
Here we can use either $\mu_I$ or  $p_I$ since 
${\rm Det}(P)$ is ${\rm SL}(2)^3$-invariant.
The formula
 simplifies considerably after we replace moments by cumulants via
(\ref{eq:takeexp})~or~(\ref{eq:momincum}):
\begin{equation}\label{eq:hyperdet3}
{\rm Det}(P)\,\,\,=\,\,\,k_{123}^{2}+4k_{12}k_{13}k_{23}. \qquad \qquad
\end{equation}
This $2 \times 2 \times 2$-hyperdeterminant is also known as the {\em tangle}, and it
appears in phylogenetics \cite{sumnerjarvis}, quantum computation \cite{miyakeQuantum} and string theory \cite{duffString}.

The next case $n=4$ is much more challenging.
According to  Huggins {\it et al.} \cite{huggins2008hyperdeterminant},
the $2{\times} 2 {\times} 2 {\times} 2$-hyperdeterminant has precisely $2,894,276$ terms,
when written as a polynomial of degree $24$ in either
probabilities $p_I$ or moments $\mu_I$. However, the expansion
of ${\rm Det}(P)$ in terms of
cumulants $k_I$ is much smaller.
The following theorem is our main result in this section.

\begin{thm} \label{thm:hyperdet}
The $2{\times} \cdots {\times} 2$-hyperdeterminant ${\rm Det}(P)$ is a polynomial
 function in the $2^n-n-1$ higher cumulants $\{k_I : |I| \geq 2 \}$.
It is homogeneous of degree $\frac{1}{2}(C_n,C_n,\ldots,C_n)$
in the $\mathbb{Z}^n$-grading given by ${\rm deg}(k_I) = \sum_{i \in I} e_i$, where $e_{i}$ is the $i$-th unit vector of $\mathbb{Z}^{n}$.
For $n = 4$, the hyperdeterminant ${\rm Det}(P)$ has precisely $13,819$ monomials in the $11$ unknowns $k_I$,
all  $\mathbb{Z}^4$-homogeneous of degree $(12,12,12,12)$, and their
total degrees range from $24$ to~$15$.
\end{thm}

 \begin{proof}
The expression of the hyperdeterminant in terms of the moments $\mu_I$
coincides with  the {\em $\mathcal{A}$-discriminant}  (cf.~\cite{gelfand1994dra})
 of the moment generating function
\begin{equation} 
\label{mxkx}
M(x) \,\,\,= \,\,\, \sum_{I \subseteq [n]}  \mu_I \prod_{i\in I} x_i  \,\,\, = \,\,\, {\rm exp}(K(x)) .
\end{equation}
Here $\mathcal{A}$ is the $(n+1) \times 2^n$ matrix whose columns are the
homogeneous coordinates of the vertices of the standard $n$-cube.
Standard  results on $\mathcal{A}$-discriminants ensure 
that ${\rm Det}(P)$ is homogeneous in the $\mathbb{Z}^{n+1}$-grading
specified by $\mathcal{A}$, so, in particular, it is homogeneous in the coarser
$\mathbb{Z}^n$-grading given by ${\rm deg}(\mu_I) = \sum_{i\in I} e_i$.
Since the degree of ${\rm Det}(P)$ in the standard $\mathbb{Z}$-grading
${\rm deg}(\mu_I) = 1$ equals $C_n$, as discussed above, we find that
${\rm Det}(P)$ is $\mathbb{Z}^n$-homogeneous of degree 
$\frac{1}{2}(C_n,C_n,\ldots,C_n)$.

The map (\ref{eq:momincum}) from moments to cumulants
respects the $\mathbb{Z}^n$-grading,
and we conclude that the expansion of ${\rm Det}(P)$ in
cumulants is  $\mathbb{Z}^n$-homogeneous of the same degree
$\frac{1}{2}(C_n,C_n,\ldots,C_n)$.
The first assertion that  ${\rm Det}(P)$ does not depend on the first order
moments $k_1,  \ldots k_n$ follows from 
Theorem~\ref{cor:invariant}.

We now come to the specific case $n=4$.
Here the proof was carried out by a computer calculation.
We first set $k_1 , k_2 , k_3 $ and $ k_4 $ to zero in the
right hand side of (\ref{mxkx})
since ${\rm Det}(P)$ does not depend on these first-order cumulants.
Our task is then to evaluate the $\mathcal{A}$-discriminant of the multilinear polynomial
\begin{eqnarray*} M(x)|_{k_1 = k_2 = k_3 = k_4 = 0} \quad = \qquad \quad 
 (k_{1234}+ k_{12}  k_{34} +k_{13}  k_{24} +k_{14}  k_{23}  )
 x_1 x_2 x_3 x_4 \\ \quad
 + \,k_{123} x_1 x_2 x_3
 + k_{124} x_1 x_2 x_4
 + k_{134} x_1 x_3 x_4
 + k_{234} x_2 x_3 x_4 \\ \quad
 + \,k_{12} x_1 x_2
 + k_{13} x_1 x_3
 + k_{14} x_1 x_4
 + k_{23} x_2 x_3
 + k_{24} x_2 x_4
 + k_{34} x_3 x_4 \,+ \,1 . \end{eqnarray*}
 This computation is done using Schl\"afli's formula
 \cite[Prop.~3]{huggins2008hyperdeterminant}. We obtained
 \begin{eqnarray*}
 256 k_{12}^6 k_{13}^5 k_{14} k_{23} k_{24}^5 k_{34}^6
-1024 k_{12}^6 k_{13}^4 k_{14}^2 k_{23}^2 k_{24}^4 k_{34}^6
+1536 k_{12}^6 k_{13}^3 k_{14}^3 k_{23}^3 k_{24}^3 k_{34}^6 + 
\\  \quad \cdots \, \hbox{many terms} \, \cdots \quad 
-\,k_{34} k_{123}^3 k_{124}^3 k_{134}^2 k_{234}^2 k_{1234}^4
+k_{123}^3 k_{124}^3 k_{134}^3 k_{234}^3 k_{1234}^3.
\end{eqnarray*}
This expansion of ${\rm Det}(P)$ has $ 13819$ terms,
all of $\mathbb{Z}^4$-degree $(12{,}12{,}12{,}12)$.
 The leading terms have total degree $24$. The last
 terms have total degree~$15$.~ \end{proof}
  
    Ideals generated by hyperdeterminants arise in various applications.
We advocate writing these in terms of cumulants. One such application, studied by
Holtz-Sturmfels \cite{holtzsturmfels} and Oeding \cite{oeding2}, concerns the
 relations among principal minors of a general symmetric $n {\times} n$-matrix $A$.
If we write $\mu_I$ for the minor with row and column indices $I \subseteq [n]$, and we 
treat the sequence $[\mu_{I}]$ as a sequence of formal moments, 
then the corresponding moment generating function  takes the special form
$$  \qquad M(x) \,\,\,=\,\,\,  {\rm det}\bigl(\,I \,+\, A  X\,\bigr) \qquad
\hbox{where} \,\,\, X \,= \, {\rm diag}(x_1,\ldots,x_n).$$
Oeding \cite{oeding2} shows that the variety of such
tables $M = [\mu_I]$ is cut out by polynomials of degree $4$.
These polynomials are obtained by acting with the group
$SL(2)^n$ on the  $2 {\times} 2 {\times} 2$-hyperdeterminants of
all subtables. We reparametrize our variety
of principal minors using the cumulant generating function:
$$ K(x) \,\,=\,\,{\rm log}\, {\rm det}(I + A X) \,\,=\,\,
{\rm trace}\, {\rm log}(I + A  X) \,\,=\,\,
{\rm trace} \bigl(\sum_{k=1}^n \frac{(-1)^{k+1}}{k} (AX)^k\bigr).
$$
The coefficients $k_I$ of the squarefree terms are sums
over all cycle monomials in $A$ that are supported on $I$. Their
algebraic relations can be computed more easily than 
those among the principal minors. We demonstrate this for  $n=4$:
$$ \begin{matrix}
K(x) & = & \sum_{I \subseteq [4]} k_I \prod_{i\in I} x_i \qquad \qquad \qquad 
\qquad \qquad \qquad \qquad \\
& = &
\sum_{i=1}^4 a_{ii} x_i
- \sum_{i<j} a_{ij}^2 x_i x_j
+ 2 \sum_{i<j<k} a_{ij} a_{ik} a_{jk} x_i x_j x_k \\  & &
- 2 (a_{12} a_{13} a_{24} a_{34} + a_{12} a_{14} a_{23}a_{34} +a_{13}a_{14}a_{23}a_{24})x_1 x_2 x_3 x_4
\end{matrix}
$$
The prime ideal of algebraic relations among the coefficients is 
found to be
% computed rapidly:
$$
\begin{small}
\!\!\!
\begin{matrix}
\bigl\langle
4 k_{12}k_{13}k_{23}+k_{123}^2,
4 k_{12}k_{14}k_{24}+k_{124}^2,
4 k_{13}k_{14}k_{34}+k_{134}^2,
4 k_{23}k_{24}k_{34}+k_{234}^2 ,\\
4 k_{12}k_{13}k_{14}k_{234}+k_{123}k_{124}k_{134}\,,\,\,
4 k_{12}k_{23}k_{24}k_{134}+k_{123}k_{124}k_{234}, \\
4 k_{13}k_{23}k_{34}k_{124}+k_{123}k_{134}k_{234}\,, \,\,
4 k_{14}k_{24}k_{34}k_{123}+k_{124}k_{134}k_{234} ,\\
2 k_{12}k_{13}k_{234}+ 2 k_{12}k_{23}k_{134} + 2 k_{13}k_{23}k_{124} + k_{123}k_{1234}, \\
2 k_{12}k_{14}k_{234}+ 2 k_{12}k_{24}k_{134}+ 2 k_{14}k_{24}k_{123} + k_{124}k_{1234}, \\
2 k_{13}k_{14}k_{234}+2 k_{13}k_{34}k_{124}+2 k_{14}k_{34}k_{123}+k_{134}k_{1234}, \\
2 k_{23}k_{24}k_{134}+2 k_{23}k_{34}k_{124}+2 k_{24}k_{34}k_{123}+k_{234}k_{1234},  \\
-2k_{12}k_{13}k_{14}k_{1234}+k_{12}k_{13}k_{124}k_{134}+k_{12}k_{14}k_{123}k_{134}+k_{13}k_{14}k_{123}k_{124},\\
-2 k_{12}k_{23}k_{24}k_{1234}+k_{12}k_{23}k_{124}k_{234}+k_{12}k_{24}k_{123}k_{234}+k_{23}k_{24}k_{123}k_{124}, \\
-2 k_{13}k_{23}k_{34}k_{1234}+k_{13}k_{23}k_{134}k_{234}+k_{13}k_{34}k_{123}k_{234}+k_{23}k_{34}k_{123}k_{134}, \\
-2 k_{14}k_{24}k_{34}k_{1234}+k_{14}k_{24}k_{134}k_{234}+k_{14}k_{34}k_{124}k_{234}+k_{24}k_{34}k_{124}k_{134},\\
k_{14}k_{123}k_{234}-k_{23}k_{124}k_{134} \,,\,\,
k_{13}k_{124}k_{234}-k_{24}k_{123}k_{134}\, , \,\,
k_{12}k_{134}k_{234}-k_{34}k_{123}k_{124}, \\
4 (k_{12}k_{13}k_{24}k_{34}{+}k_{12}k_{14}k_{23}k_{34}{+} k_{13}k_{14}k_{23}k_{24}) \qquad \\
\qquad {-}2 (k_{14}k_{123}k_{234}{+} k_{24}k_{123}k_{134} {+}  k_{34}k_{123}k_{124}) {-} k_{1234}^2
\bigr\rangle
\end{matrix}
\end{small}
$$
These twenty polynomial correspond to the hyperdeterminantal 
relations in \cite[Thm.~8]{holtzsturmfels}. They furnish
a compact encoding  of this codimension $5$ variety.

\section{Invariance and Independence}\label{sec:unipotent}

The algebraic relations in Section 2
did not involve any of the order one cumulants $k_1,\ldots,k_n$
and they were homogeneous with respect to the
$\mathbb{Z}^n$-grading given by ${\rm deg}(k_I) = \sum_{i\in I} e_i$.
In this section we argue that these properties  hold for  all
statistically meaningful varieties in the space of
$2 {\times} \cdots {\times} 2$-tables.

To compute moments in (\ref{eq:mominprob}) we used the convention that the (formal) random vector $X=(X_{1},\ldots, X_{n})$ has values in $\{0,1\}^{n}$.  Other authors prefer the choice $\{-1,1\}^{n}$, and 
this leads to rather different formulas for the moments (see \cite[Equation (2.1)]{bruno1999probability}).
A meaningful statistical model will not depend on such choices. Hence we are only interested
in cumulant varieties that do not depend on such choices.

Suppose we replace each of our random variable $X_i$ by
a new variable $X_i'$ which takes values $a_i$ and $b_i$ instead of $1$  and $0$.
If the probability distribution is
the same on both state spaces,
then the cumulants are transformed via 
\begin{equation}\label{eq:cummultipl} \qquad
k_{I}'\,\,\, = \,\,\, k_I \cdot \prod_{i\in I}(a_i-b_i) \qquad\mbox{ for all } I\subseteq [n] 
\hbox{ and } |I|\geq 2
\end{equation}
and $k_{i}'=(a_{i}-b_{i})k_{i}+b_{i}$ for $i=1,\ldots, n$. 
%See e.g.~\cite{mccullagh1987tms}.
This result is purely algebraic and the above remains true if we replace probability distributions with any complex distributions. In geometric language, changing the values of the binary
variables $X_i$ corresponds to a natural action
of the $n$-dimensional torus
$(\mathbb{C}^*)^n$ with coordinates $a_i-b_i$ on
the space  $\mathbb{C}^{2^n-n-1}$ whose
coordinates are the higher cumulants $k_I$, $|I| \geq 2$.
This action is compatible with the $\mathbb{Z}^n$-grading:

\begin{thm} \label{thm:invariant}
A subvariety of $\mathbb{C}^{2^n-1}$ is
invariant under changing values of the $X_i$
if and only it is defined by  $\mathbb{Z}^n$-homogeneous polynomials
in $k_I$ with~$|I| \geq 2$.
\end{thm}

\begin{proof}
Let $V$ be a subvariety in the space $\mathbb{C}^{2^n-1}$ 
whose coordinates are all the cumulants. Suppose that $V$
is invariant under replacing the values $(0,1)$ of $X_i$
by any $(b_i,a_i)$. If the new values satisfy $a_i= b_i+1$ then
the higher cumulants $k_I$, $|I| \geq 2$, remain unchanged
but the vector $(k_1,\ldots,k_n)$ is shifted to
$(k_1 + b_1, \ldots, k_n+b_n)$.
Hence the ideal $I_V$ of $V$ is generated by
polynomials that do not depend on linear cumulants $k_1,\ldots,k_n$.
By fixing $b_i = 0$ and moving $a_i$, we see that
$V$ is invariant under the torus action (\ref{eq:cummultipl}). Hence
its ideal $I_V$ is $\mathbb{Z}^n$-homogeneous, and this
proves the only-if direction. The if-direction holds by essentially
the same argument.
\end{proof}

The group $SL(2)^n$ acts on the
tensor space $\mathbb{C}^{2 \times 2 \times \cdots \times 2}$
and many important varieties are
invariant under this action. In particular, they are
invariant under $U(2)^n$ where $U(2)$ is the unipotent group
of $2\times 2$-matrices of the form
$$
\left[\begin{array}{cc}
1 & \lambda\\
0 & 1
\end{array}
\right]\qquad \mbox{for }\lambda\in \mathbb{C}.
$$ 
The invariance property of Theorem \ref{thm:invariant}
reflects precisely the $SL(2)^n$-invariance.

\begin{cor} \label{cor:invariant}
Let $V$ be a subvariety of the affine open subset $\{\mu_\emptyset = 1\}$ in 
the projective space $\mathbb{P}(\mathbb{C}^{2 \times 2 \times \cdots \times 2})$
and let $\overline{V}$ denote its closure in that projective space.
If  $\overline{V}$ is invariant
under the action of $SL(2)^n$ then the ideal $I_V$ that defines $V$ is generated by
$\mathbb{Z}^n$-homogeneous polynomials in the $k_I$ with $|I| \geq 2$.
\end{cor}

\begin{proof}
The unipotent group $U(2)^n$ acts on the moment generating function
 via $M(x)\mapsto M(x)\prod_{i=1}^{n}(1+\lambda_{i}x_{i})$. 
 Modulo the ideal  $\langle x_{1}^{2},\ldots, x_{n}^{2}\rangle$
 we have
 $$
M(x)\prod_{i=1}^{n}(1+\lambda_{i}x_{i})
\quad = \quad M (x) \exp\bigl(\sum_{i=1}^{n}\lambda_{i}x_{i}\bigr)
\quad = \quad \exp \bigl (K(x)+\sum_{i=1}^{n}\lambda_{i}x_{i} \bigr).
$$
This means that $U(2)^{n}$ acts on the space of cumulants by shifting the first order cumulants.
We conclude that the prime ideal of $V$ 
is generated by polynomials in the cumulants $k_I$ with $|I| \geq 2$.
Since $V$ is also invariant under tuples of $2 \times 2$-diagonal matrices in $SL(2)^n$,
these ideal generators can be chosen to be $\mathbb{Z}^n$-homogeneous.
\end{proof}

Hyperdeterminants and their ideals  in Section 2
are $SL(2)^n$-invariant and hence 
expressible
%representable
 by
$\mathbb{Z}^n$-homogeneous polynomials in higher~cumulants.

\begin{exmp}
The converse does not hold in Corollary \ref{cor:invariant}.
Fix  $n = 3$, let $ \rho \in \mathbb{C} \backslash \{4\}$, and consider the hypersurface  
in $\{ \mu_\emptyset = 1\} \subset \mathbb{P}(\mathbb{C}^{2 \times 2 \times 2})$
 defined~by 
$$ k_{123}^2 \,+\, \rho \cdot k_{12} k_{13} k_{23} \,\,\, = \,\,\,  0 .$$
This equation has degree six when written in the (homogenized) moments:
$$ {\rm Det}(P) m_\emptyset^2 \,+\,(\rho-4)
(m_\emptyset m_{12} -m_1 m_2)
(m_\emptyset m_{13} -m_1 m_3)
(m_\emptyset m_{23} -m_2 m_3)
$$
This defines a sextic hypersurface in $\mathbb{P}(\mathbb{C}^{2 \times 2 \times 2})$ 
that is $U(2)^3$-invariant
but {\em not $SL(2)^3$-invariant}.
The formula in probabilities is even less invariant:
\begin{eqnarray*}
{\rm Det}(P) \cdot
(p_{\emptyset}+p_1+p_2+p_3+p_{12}+p_{13}+p_{23}+p_{123})^2 
\qquad \qquad \\
\, + \,
(\rho -4) 
(p_\emptyset  p_{23}+p_\emptyset p_{123} + 
p_1 p_{23} + p_1 p_{123} - p_2 p_3 - p_2 p_{13} - p_3 p_{12}  - p_{12} p_{13})\\
 \cdot
( p_\emptyset p_{13} + p_\emptyset p_{123} - p_1 p_3 -
p_1 p_{23} + p_2 p_{13} + p_2 p_{123} - p_3 p_{12} - p_{12} p_{23}) 
\\ 
\cdot
(p_\emptyset p_{12} + p_\emptyset p_{123} - p_1 p_2 - p_1 p_{23}
-p_2 p_{13} + p_3 p_{12}
+p_3 p_{123} - p_{13} p_{23})
\end{eqnarray*}
Of course, for $\rho = 4$, this is
the hyperdeterminantal quartic $\{ {\rm Det}(P) = 0 \}$. \qed
\end{exmp}

The most basic statistical model for $n$ binary random variables  $X_i$ is 
the model of {\em complete independence}, denoted $X_{1}\indep X_{2}\indep \ldots\indep X_{n}$,
which is the Segre variety
$(\mathbb{P}^1)^n \subset \mathbb{P}^{2^n-1}$.
In terms of moments, this is parametrized by
$\,M(x) = \prod_{i=1}^n (1+ \mu_i x_i)  $. In terms of cumulants, we obtain
$\,K(x) =  \sum_{i=1}^n {\rm log}(1+\mu_ix_i) =  \sum_{i=1}^n k_i x_i$.
In probability coordinates, the Segre variety is defined by certain
$2 {\times} 2$-determinants $p_I p_J - p_K p_L$
but we see that this simplifies when we use cumulants as coordinates:

\begin{rem} \label{rem:segre}
The  Segre variety is defined by
$k_I = 0 $ for all $|I| \geq 2$.
\end{rem}

The Segre variety is the intersection of the independence models
$A \indep B$ where $A|B$ runs over all partitions of the set $[n]$.
The equations for $A \indep B$
are $k_I = 0$ for all $I$ with $A \cap I \not= \emptyset $
and $B \cap I \not= \emptyset$. 
The model $A \indep B$ also makes
sense when $A \cup B$ is a proper subset of $[n]$,
with equations  as follows.

\begin{prop}
If $A$ and $B$ are disjoint subsets of $[n]$ then the
independence model $A \indep B$
is defined by
$k_I=0$ where $I \subseteq A \cup B$, $A \cap I \not= \emptyset $
and $B \cap I \not= \emptyset$. 
\end{prop}

\begin{proof}
The independence model $A \indep B$ has the moment parametrization 
$$ M(x) \quad = \quad M_1(x_i: i \in A) \cdot M_2(x_j : j \in B) \, \,\, + \!\!\! \sum_{l \in [n] \backslash (A \cup B)} 
\!\!\!\!\! x_l \cdot N_l (x) . $$
By taking the logarithm, we find that
$$ K(x) \,= \, {\rm log}(M_1(x_i:i \in A)) \, + \,{\rm log}(M_2(x_j: j \in B))
\,\, {\rm mod} \,\,\langle x_l: l \in  [n] \backslash (A {\cup B)} \rangle.
$$
This form is equivalent to the asserted vanishing condition on cumulants.
\end{proof}

The symmetry group of the $n$-cube is the semidirect product of the symmetric group $S_n$,
which permutes $[n]$, and the abelian group $\mathbb{Z}_2^n$,
which swaps $0$s and $1$s.
This gives rise to an action on 
 $\mathbb{C}^{2 {\times} 2 {\times} \cdots {\times} 2}$.
We identify elements $\rho\in \mathbb{Z}_{2}^{n}$ with subsets $J\subseteq [n]$.
The action on coordinates $p_I$ is as follows:
for $\sigma\in S_{n}$ we have $\sigma (p_{I})=p_{\sigma(I)}$,
and for $ J\subseteq [n]$ we have
 $\rho_{J} (p_{I})=p_{I\Delta J}$, where $I\Delta J=(I \backslash J)\cup 
 (J \backslash I)$. This being an action ensures
\begin{equation}
\label{actionensures}
 \rho_{J}(p_{I}) \,=\,\rho_{j_{r}}\circ\cdots \circ\rho_{j_{1}}(p_{I}) \quad
\hbox{for all} \quad J=\{j_{1},\ldots,j_{r}\}\,\subseteq \,[n].
\end{equation}

We extend this action from probability coordinates
to any of their polynomials in a natural way. In this way we extend this action to cumulant coordinates. This action is simple for permutations
 $\sigma\in S_{n}$: we have  $\sigma (k_{I})=k_{\sigma(I)}$. The action of 
 the group $\mathbb{Z}_{2}^{n}$ is more subtle,
 and it can be characterized by the following corollary.
 That result will help us in  Section \ref{sec:semialg} to get a more compact semialgebraic 
description of the space of cumulants, by
taking advantage of  the symmetries in our problem.

\begin{cor}\label{cor:sym-dif}
Consider the cumulants $k_I$
as polynomials in probabilities $p_I$, via
 (\ref{eq:mominprob}) and (\ref{eq:cum1}).
 For $I,J\subseteq [n]$ with $|I|\geq 2$,
 the element $\rho_J  \in \mathbb{Z}_2^n$ satisfies
\begin{equation} \label{eq:swaplabels}
\rho_J (k_{I}) \,\,\, = \,\,\,
\left\{\begin{array}{ll}
-k_{I} & \mbox{ if } \, |J\cap I| \mbox{ is odd},\\
\phantom{-} k_{I} & \ \mbox{otherwise.}
\end{array}\right.
\end{equation}
Furthermore, for each $i=1,\ldots,n$, we have
$$
\rho_J(k_{i}) \,\,\, = \,\,\, \left\{\begin{array}{ll}
1-k_{i} & \mbox{if }\, i\in J,\\
\,\,\, \, k_{i} & \mbox{otherwise.}\end{array}\right.
$$
\end{cor}

\begin{proof}
By (\ref{actionensures}) it suffices to show $\rho_i(k_{I})=-k_{I}$ if
$i\in I$ and $\rho_i(k_{I})=k_{I}$ if $i\notin I$. 
Formula (\ref{eq:swaplabels}) follows from (\ref{eq:cummultipl})
by taking $a_i = 0 $ and $b_i = 1$, i.e.~we swap the states
of the $i$th variable $X_i$, and similarly for first-order cumulants.
\end{proof}

\section{Tangents and secants of the Segre variety}

In Remark \ref{rem:segre} we saw that the Segre variety  $\,(\mathbb{P}^1)^n
\subset \mathbb{P}(\mathbb{C}^{2 \times \cdots \times 2}) $
collapses to a single point  in the space $\mathbb{C}^{2^n-n-1}$
of higher cumulants. This raises the question what the representation in
the $k_I$ with $|I|\geq 2$ looks like for
varieties naturally associated to $(\mathbb{P}^1)^n$, such as  its secant and tangential varieties.
 We here examine the first tangential variety and the first secant variety:
 $$
\begin{matrix}
{\rm Tan}((\mathbb{P}^1)^n)&= & \hbox{closure of} & \{\,x\in \mathbb{P}^{2^{n}-1}\,|\, \,x \mbox{ lies on a line tangent to } (\mathbb{P}^1)^n\}, \\
{\rm Sec}((\mathbb{P}^1)^n)& = & \hbox{closure of} & \{\,x\in \mathbb{P}^{2^{n}-1} \,|\, \, x \mbox{ lies on a secant line of } (\mathbb{P}^1)^n\}.
\end{matrix}
$$
Our next result reveals that the tangential variety is toric
in the cumulants.

\begin{thm} \label{thm:tangtoric}
The image of the tangential variety ${\rm Tan}((\mathbb{P}^1)^n)$ in the
 space of higher cumulants $\,\mathbb{C}^{2^n-n-1}\,$ 
 is isomorphic to the $n$-dimensional affine toric variety parametrized by all
square-free monomials of degree $\geq 2$.
\end{thm}

\begin{proof}
In the tensor notation of \cite{oeding1},  ${\rm Tan}((\mathbb{P}^1)^n)$ has the parametrization
$$ M \quad = \quad\frac{1}{n}
\sum_{i=1}^{n} a^{{(1)}}\otimes\cdots\otimes a^{{(i-1)}}\otimes b^{{(i)}} \otimes a^{{(i+1)}}\otimes \cdots \otimes a^{{(n)}},$$
where $a^{(i)} = (1,a_i)$ and  $b^{(i)} = (1,b_i)$
are vectors representing points in the distinguished affine open subset of $\mathbb{P}^1$. The formula above translates into the following parametrization of moment generating functions:
\begin{equation}\label{eq:tangmoment}
M(x) \quad = \quad \prod_{j=1}^n (1+ a_j x_j) \cdot \biggl( \sum_{i=1}^n \frac{1 + b_i x_i}{1+ a_i x_i} \biggr) .
\end{equation}
We compute the logarithm
of the series $M(x)$ modulo $\langle x_1^2,x_2^2,\ldots,x_n^2 \rangle$.
Disregarding $\mathbb{R}$-linear combinations of $x_1,\ldots,x_n$, and
 setting $s_i = (a_i-b_i)/n $, 
 %we find
 \begin{equation}
 \label{eq:nontrivial}
 K(x) \,\,= \, \,{\rm log} \bigl( \sum_{i=1}^n \frac{1 + b_i x_i}{1+ a_i x_i} \bigr)
\quad = \quad \sum_{ |I| \geq 2} (-1)^{| I |-1} (|I|-1) ! \prod_{i \in I} s_i x_i ,
\end{equation}
The identity on the right can be proved directly by manipulating
generating functions. An alternative and more detailed proof will be given in
Example \ref{ex:tangentialrevisited}.
 We now conclude that
\begin{equation}
\label{eq:toricpara}
\qquad   k_I \,\,\, = \,\,\,  (-1)^{| I |-1} (|I|-1) ! \cdot  \prod_{i \in I} s_i \qquad \hbox{for} \quad | I | \,\geq \,2.
 \end{equation}
 This monomial parametrization shows that ${\rm Tan}((\mathbb{P}^1)^n)$ is  toric
   in cumulants. It is isomorphic to the toric variety with
 parametrization $\,  k_I \,\mapsto \, \prod_{i \in I} s_i  $.
\end{proof}

We easily find the cumulant ideal of ${\rm Tan}((\mathbb{P}^1)^n)$,
by computing a Markov basis for the toric ideal of relations among
squarefree polynomials of degree $\geq 2$.
We then  rescale to adjust to the signs and factorials appearing in~(\ref{eq:toricpara}).

\begin{exmp}
Let $n=5$. Then the toric ideal ${\rm Tan}((\mathbb{P}^1)^n)$ 
is minimally generated by $120$ binomials
in the $26$ cumulant coordinates.
Among these generators, $75$ are quadrics
and $45$ are cubics. The quadrics include binomials such as
$k_{12} k_{34}-k_{14} k_{23}$, 
$k_{123} k_{45}-k_{12} k_{345}$,
$k_{123} k_{345} - k_{135} k_{234}$,
$k_{1234}+6 k_{14} k_{23}$, and
$ k_{12345}+12 k_{12} k_{345}$.
The cubics include binomials such as
$\,k_{123}^2+4 k_{12} k_{13} k_{23} \,$ and
$\,k_{123} k_{124}+4 k_{12} k_{14} k_{23}$. \qed
\end{exmp}

We now come to the {\em secant variety} $\,{\rm Sec}((\mathbb{P}^1)^n)$. This
is not a toric variety in cumulants. For example, for $n=4$ it has
the following parametrization:
$$
\begin{matrix}
M & = &  (1-t)\,A \otimes B \otimes C \otimes D \,\,\,+\,\,\,
t\,E \otimes F \otimes G \otimes H, \\
M(x) & = &
 (1-t)(1{+}ax_1)(1{+}bx_2)(1{+}c x_3)(1{+}d x_4) \\ & &  \,\,+\,t \,
 (1{+}ex_1)(1{+}fx_2)(1{+}g x_3)(1{+}h x_4).
 \end{matrix}
$$
 Here $A = (1,a) ,\ldots, H = (1,h)$, and $t$ is  a complex mixing parameter.
  
The image of ${\rm Sec}((\mathbb{P}^1)^n)$ in the $11$-dimensional 
space of higher cumulants is a $5$-dimensional
affine variety that is not toric. Its ideal is
generated by $16$ polynomials in $k_{12}, k_{13}, 
\ldots,k_{1234}$. These are 
the ten binomial quadrics
\begin{equation}
\label{eq:GENS1}
\begin{matrix}
k_{12} k_{34} - k_{14} k_{23}, k_{13} k_{24} - k_{14} k_{23}, \,
k_{12} k_{134} - k_{14} k_{123} , k_{13} k_{124} - k_{14} k_{123},   \\
k_{12} k_{234} - k_{24} k_{123} \,,\,\, k_{23} k_{124} - k_{24} k_{123} \, ,   \,\,
k_{13} k_{234} - k_{34} k_{123}, \\ k_{23} k_{134} - k_{34} k_{123} \,, \,\,
k_{14} k_{234} - k_{34} k_{124} \,, \,\, k_{24} k_{134} - k_{34} k_{124}, 
\end{matrix}
\end{equation}
and the six non-binomial cubics
\begin{equation}
\label{eq:GENS2}
\begin{matrix}
k_{12}L - (k_{123} k_{124} + 4 k_{12} k_{14} k_{23}) , &
k_{13}L - (k_{123} k_{134} + 4 k_{13} k_{14} k_{23}), \\
k_{14}L  - (k_{124} k_{134} + 4 k_{14} k_{14} k_{23}), &
k_{23}L - (k_{123} k_{234} + 4 k_{23} k_{14} k_{23}), \\
k_{24}L - (k_{124} k_{234} + 4 k_{24} k_{14} k_{23}), &
k_{34}L  - (k_{134} k_{234} + 4 k_{34} k_{14} k_{23}).
\end{matrix}
\end{equation}
Here $\,L = k_{1234} + 6 k_{14} k_{23} \,$ is one of the
toric relations on (\ref{eq:toricpara}).
Indeed,  the tangential variety ${\rm Tan}((\mathbb{P}^1)^4)$ is a hypersurface
in the secant variety ${\rm Sec}((\mathbb{P}^1)^4)$. Its toric ideal in
cumulants has $21$ minimal generators,
namely, the ten quadrics in (\ref{eq:GENS1}),
the quadric $L$, the six parenthesized cubics in (\ref{eq:GENS2}),
and the four  hyperdeterminants
$$
    k_{234}^2 + 4 k_{23} k_{24} k_{34} \, , \,\, 
    k_{134}^2 + 4 k_{13} k_{14} k_{34} \, , \,\,
    k_{124}^2 + 4 k_{12} k_{14} k_{24} \, , \,\,
    k_{123}^2 + 4 k_{12} k_{13} k_{23}.
$$
These various equations in cumulants can now be
translated back into probability coordinates,
using the substitutions (\ref{eq:mominprob}) and (\ref{eq:cum1}).
After homogenizing and saturating with $\mu_\emptyset$, we recover
the $32$-dimensional space of $3 {\times} 3$-minors
of flattenings for ${\rm Sec}((\mathbb{P}^1)^4)$,
as in \cite{raicu},
and the $53$ ideal generators for ${\rm Tan}((\mathbb{P}^1)^4)$,
namely, the $32$ cubics, the $20$ hyperdeterminantal quartics,
and  the special quadric as in \cite[\S 3.2]{oeding1}.

For any $n \geq 4$,  the secant variety ${\rm Sec}((\mathbb{P}^1)^n)$ is a  curve over the toric variety
 ${\rm Tan}((\mathbb{P}^1)^n)$.
Using cumulant coordinates,    it has the parametrization 
 \begin{equation}
 \label{eq:secpara}
 k_I \quad = \quad \kappa_{|I|} (t) \cdot \prod_{i \in I } b_i , 
 \end{equation}
  where $b_{i}$ are complex parameters and $\kappa_\nu(t)$ is a certain univariate polynomial of degree $\nu$;
  see (\ref{eq:necklaces}).
  For example, 
 \begin{equation}
 \label{eq:kappa}
 \begin{matrix} 
  % \kappa_1(t) &=&   t ,\\
  \kappa_2(t) &=& -t^2+t  \\
\kappa_3(t) &=&   2t^3-3t^2+t , \\
 \kappa_4(t) &=& -6 t^4+12t^3-7t^2+t.
 \end{matrix}
\end{equation}
The leading coefficient of $\kappa_\nu(t)$ equals $(-1)^{\nu-1} (\nu-1)!$
in the parametrization  (\ref{eq:toricpara}) of the tangential variety.
Using (\ref{eq:kappa}), we can now   recover the equations
 (\ref{eq:GENS1}) and (\ref{eq:GENS2})
 of the secant variety by
implicitizing (\ref{eq:secpara}) for $n=4$.
The derivation of  (\ref{eq:secpara})
and the polynomials $\kappa_\nu(t)$
will be explained in Example~\ref{ex:secantrevisited}.

\section{Hidden subset models}

We now introduce a highly overparametrized algebraic statistical model for a vector $X$
of $n$ binary random variables.
It is called the {\em complete hidden subset model}. Here is a generative description of this model.
A subset $I$ of $[n]$ (or alternatively a binary vector) is to be chosen at random.  For each element $i \in [n]$ we need to 
decide whether $i$ is in $I$ or not. This is done as follows.
First, a hidden subset $J$ is chosen with some probability $t_J$.
Then we select $i$ for $I$ with probability $a_i^{(0)}$ if $i \not\in J$, and we select $i$ 
for $I$ with probability $a_i^{(1)}$ if $i \in J$. The conditional probabilities
$\,a_i^{(0)} = {\rm Prob}( i\in I\, |\, i \not\in J )\,$ and
$\,a_i^{(1)} = {\rm Prob}( i\in I \,|\, i \in J )\,$
are unrelated parameters that govern this process.

The distributions  in this model are parametrized as follows:
$$
p_{I} \,\,\, =\,\,\, \sum_{J\subseteq [n]} t_{J} \prod_{i\in I^{c}\cap J^{c}}(1-a_{i}^{(0)})\prod_{i\in I^{c}\cap J}(1-a_{i}^{(1)})\prod_{i\in I\cap J^{c}}a_{i}^{(0)}\prod_{i\in I\cap J}a_{i}^{(1)},
$$
where $I^{c}$ denotes the complement of $I\subseteq [n]$.
The corresponding moment generating function has the parametrization
   \begin{equation}\label{eq:momsplits}
M(x)\quad = \quad
\sum_{J\subseteq [n]} t_{J} \cdot \prod_{i\in J^{c}}(1+a_{i}^{(0)}x_{i}) \prod_{i\in J}(1+a_{i}^{(1)}x_{i}).
\end{equation}
 
 The model has two parameters $a_{j}^{(0)}$ and $a_{j}^{(1)}$,
with values between $0$ and $1$, for each $j \in [n]$.
Further, it has $2^{n}$ mixing parameters $t_{I}$,
one for each subset $I\subseteq [n]$.
These parameters are non-negative and they sum to $1$,
so the tables $T = (t_I)_{I \subseteq [n]}$ is also a distribution. We write $k_I^{(t)}$ for the cumulants obtained from the table $T$.

\smallskip

Our main result in this section is the following  intriguing theorem.

\begin{thm}\label{thm:bigCSImodel}
The complete hidden subset model is parametrized in terms of cumulants by 
$\,k_{i}=a_{i}^{(0)}+b_{i} \cdot k^{(t)}_{i}$,
where $\,b_{i}=a^{(1)}_{i}-a^{(0)}_{i}\,$ for $i=1,\ldots,n$, and
\begin{equation}
\label{eq:remarkable}
k_{I} \quad = \quad k_{I}^{(t)} \cdot \prod_{i\in I}b_{i}\qquad\mbox{for }\,\, |I|\geq 2.
\end{equation}
\end{thm}

\begin{proof}
We introduce a homogeneous probability generating function as
$$
P_{{\rm hom}}(y^{(0)},y^{(1)})=\sum_{J\subseteq [n]} p_{J}\prod_{i\in J^{c}}y_{i}^{(0)}
\prod_{i\in J}y_{i}^{(1)},
$$
so that $P(x)=P_{{\rm hom}}(\mathbf{1},x)$.
Then the moment generating of $X$ in (\ref{eq:momsplits}) can be dually treated as 
the homogeneous version of the probability generating function of $Y$. Namely, for fixed $a_{i}^{(0)}$ and $a_{i}^{(1)}$, we write 
(\ref{eq:momsplits})~as $M(x)= P^{(t)}_{{\rm hom}}(y^{(0)},y^{(1)})$, 
where $y_{i}^{(0)}=1+a_{i}^{(0)}x_{i}$ and $y_{i}^{(1)}=1+a_{i}^{(1)}x_{i}$. From the homogeneous  generating function $P^{(t)}_{{\rm hom}}$ we can obtain the homogeneous moment generating function $M^{(t)}_{{\rm hom}}$ similarly as in the first equation in (\ref{eq:mgf}). Thus setting $z_{i}=y^{(1)}_{i}-y^{(0)}_{i}$, we have
\begin{equation}\label{eq:hompgf}
P^{(t)}_{{\rm hom}}(y^{(0)},y^{(1)})\quad=\quad P^{(t)}_{{\rm hom}}(y^{(0)},z+y^{(0)})\quad = \quad M^{(t)}_{{\rm hom}}(y^{(0)},z).
\end{equation}
From this we find 
that $M(x)=P^{(t)}_{{\rm hom}}(y^{(0)},y^{(1)})$ is equal to 
$$
M^{(t)}_{{\rm hom}}(y^{(0)},z)\quad =\quad\sum_{J\subseteq [n]} \mu_{J}^{(t)}\prod_{i\in J}z_{i}\prod_{i\in J^{c}}y_{i}^{(0)}.
$$
Since $z_{i}=y^{(1)}_{i}-y^{(0)}_{i}=b_{i}x_{i}$ and $y_{i}^{(0)}=1+a_{i}^{(0)}x_{i}$, then $M^{(t)}_{{\rm hom}}(y^{(0)},z)$ can be rewritten as a function of $x_{1},\ldots,x_{n}$ only: 
$$
\begin{array}{rcl}
M^{(t)}_{{\rm hom}}(y^{(0)},z)&=&\sum_{J\subseteq [n]} \mu_{J}^{(t)}\prod_{i\in J}b_{i}x_{i}\prod_{i\in J^{c}}(1+a^{(0)}_{i}x_{i})=\\
\quad&=&\quad {M}^{(t)}(b_{1}x_{1},\ldots,b_{n}x_{n})\prod_{i=1}^{n}(1+a^{(0)}_{i}x_{i}).
\end{array}
$$
The last equality follows holds modulo the ideal $\langle x_{1}^{2},\ldots,x_{n}^{2}\rangle$. This implies $K(x)={K}^{(t)}(b_{1}x_{1},\ldots,b_{n}x_{n})+\sum_{i=1}^{n} a_{i}^{(0)}x_{i}$ and hence $k_{i}=a_{i}^{(0)}+b_{i}k^{(t)}_{i}$ for $i=1,\ldots, n$, and $k_{I}=k_{I}^{(t)}\prod_{i\in I} b_{i}$ for every $I\subseteq [n]$ with $|I|\geq 2$.
\end{proof}

%\begin{rem}A statistician may note that
%(\ref{eq:remarkable}) can be derived from the formula in \cite{brillinger1969calculation} for
 % cumulants in terms of the conditional cumulants.  If the $X_i$ are independent given 
 % the mixing variable $Y$ then all the joint conditional cumulants vanish, except for the conditional means. 
%  The result follows by checking that $\mathbb{E}(X_{i}|Y)=b_{i}(Y_{i}-\mathbb{E} Y_{i})$.
%\end{rem}

A {\em hidden subset model} is any submodel obtained from (\ref{eq:momsplits}) by
setting some of the mixing parameters $t_{I}$ to zero.
Thus a hidden subset model for $n$ binary variables
is specified by a collection $\{I_1,\ldots,I_k\}$ of subsets of $[n]$.
These subsets indicate those mixing parameters $t_{I_1}, \ldots, t_{I_k}$
that are not zero. We next show that the two varieties in  Section~4
arise as special cases of this.

\begin{exmp}
\label{ex:tangentialrevisited}
The hidden subset model  $\bigl\{ \{1\}, \ldots, \{n\} \bigr\}$
is given by
$$ M(x) \quad = \quad \prod_{j=1}^n (1+ a_j^{(0)} x_j) \cdot \biggl( \sum_{i=1}^n t_{i} 
\frac{1 + a_i^{(1)} x_i}{1+ a_i^{(0)} x_i} \biggr). $$
This equals the tangential variety 
${\rm Tan}((\mathbb{P}^1)^n)$ as in (\ref{eq:tangmoment}) but now with
toric parameters  
$s_i = t_i  (a_i^{(1)}- a_i^{(0)}) $.
%$s_i = t_i  (a_i^{(0)}- a_i^{(1)})/n $. % POITR, PLEASE CHECK
We compute the cumulants $k_I^{(t)}$ of the mixing distribution
$T$ using   (\ref{eq:cum1}). The moments of $T$ satisfy
$\mu_{i}^{(t)}=t_{i}$ and $\mu_{I}^{(t)}=0$ for $|I|\geq 2$.
This means that the sum in (\ref{eq:cum1}) has
 only one non-zero term, the one
corresponding to the partition $\pi$ of $I$ into singleton blocks.
Now,  (\ref{eq:cum1}) reads
$$
k_{I}^{(t)}\,\,\,= \,\,\, (-1)^{|I|-1}(|I|-1)! \prod_{i \in I} t_i
$$
We have shown that the formula (\ref{eq:remarkable}) specializes
to (\ref{eq:toricpara}) for  this model. \qed
\end{exmp}

\begin{exmp}
\label{ex:secantrevisited}
The hidden subset model  $\bigl\{\emptyset, [n] \bigr\}$
is the mixture of two independent random variables,
so it coincides with the secant variety 
${\rm Sec}((\mathbb{P}^1)^n)$.
The mixing distribution $T$ has one free mixing parameter $t$, where
 $  t_{1\cdots n} = t $ and $t_\emptyset = 1-t$. The moments of $T$
are $\mu_\emptyset = 1$ and $\mu_B^{(t)} = t$ for $|B| \geq 1$.
The formula (\ref{eq:cum1})  implies
$$ k_I^{(t)} \,\,\,= \,\,\, \sum_{i=1}^{|I|} \, (-1)^{i-1} (i-1)! \cdot \beta_{i,| I |} \cdot t^i  $$
where $\beta_{i,I}$ is the number of set partitions
of $I$ into $i$ blocks.
This univariate polynomial depends only on the cardinality $\nu = |I|$.
We can also write it~as 
\begin{equation}
\label{eq:necklaces}
\kappa_\nu (t) \,\, = \,\,\sum_{i=1}^\nu \, (-1)^{i-1}  \cdot \gamma_{i, \nu} \cdot t^i 
\end{equation}
where $\gamma_{i,I}$ is the number of cyclically ordered set partitions
of a $\nu$-set into $i$ blocks.
Such partitions are known as {\em necklaces} in enumerative combinatorics. \qed
\end{exmp}

\begin{exmp}[Binary Hidden Markov Model]
The complete hidden subset model includes all {\em hidden Markov models} (HMM)
where both  hidden and  observed states are binary.
These models are widely used in computational biology \cite[\S 1.4.3 and \S 11]{ascb}.
We can treat the mixing variable with distribution $T$ as a hidden binary process $Y=(Y_{1},\ldots,Y_{n})$. 
The parameters $a_{i}^{(0)} $ and $a_{i}^{(1)}$ determine the conditional distribution of $X_{i}$ given the hidden process, and this observed distribution depends on $Y$ only through the value of $Y_{i}$. In this context the parameters $b_{i}=a_{i}^{(1)}-a_{i}^{(0)}$ are the linear regression coefficients of $X_{i}$ with respect to $Y_{i}$.
For an HMM, the hidden distribution $T$ follows, in addition, a homogeneous Markov chain
\cite[\S 1.4.2]{ascb}. Thus, if $k_I^{(t)}$ are the cumulants of the 
homogeneous Markov chain,
then (\ref{eq:remarkable}) gives a parametrization of the binary HMM. It would be interesting
to revisit the recent work of Sch\"onhuth \cite{schoenhuth} from this perspective.
We expect his prime ideals $I_{3,n}$ in \cite[\S 7.3]{schoenhuth} to
 have a nice representation  in terms of cumulants. \qed
\end{exmp}

The set of all hidden subset models, for fixed $n$, forms a poset
whose elements $\mathcal{M}_A$ are indexed by the
$2^{2^n}$ subsets $A$ of $2^{[n]}$. The model $\mathcal{M}_A$
is obtained from the complete hidden subset model by
setting $t_{I}=0$ for all $I\subseteq [n]$  not in $A$.
Of course, different labels  $A$ and $B$ can lead to isomorphic 
hidden subset models
$\mathcal{M}_{A}$ and $\mathcal{M}_{B}$. Clearly,
this happens if $B$ is obtained from $A$ by a permutation of $[n]$. 
But, in fact, the full symmetry cube of the $n$-cube acts
on the hidden subset models:

\begin{prop} \label{prop:symmetry}
Let $A,B\subseteq 2^{[n]}$ and assume that
$B$ is equal to $J \Delta A =\{ I\Delta \alpha : \alpha \in A\}$
for some subset  $J\subseteq [n]$.
Then $\mathcal{M}_{A}$ is isomorphic to $\mathcal{M}_{B}$.
\end{prop}
\begin{proof}
By Theorem \ref{thm:bigCSImodel}, the two models are parametrized by $k_{I}=k_{I}^{(t)}\prod_{i\in I} b_{i}$.
 The cumulants $k_{I}^{(t)}$ depend on $A$ and $B$. By Corollary \ref{cor:sym-dif},
if $B=J\Delta A$ for some $J$, then the respective cumulants
 $k_{I}^{(t)}$ for $A$ and $B$ agree up to sign.
 \end{proof}

Each hidden subset model $\mathcal{M}_{A}$ can be 
identified with a $0/1$-polytope $P_{A}$.
%Probability assignments for the mixing variable $Y$
%determine points of $P_{A}$ in an obvious way. 
By Proposition \ref{prop:symmetry}, if $P_{A}$ and $P_{B}$ are $0/1$ equivalent then $\mathcal{M}_{A}$ and $\mathcal{M}_{B}$ are isomorphic. 
We say that the model $\mathcal{M}_{A}$ is {\em non-degenerate} if 
the polytope $P_A$ is not contained in any hyperplane $x_{i}=0$ or $x_{i}=1$ for $i=1,\ldots,n$.
If this happens then the random variable $X_{i}$ is independent of all  other variables. Geometrically this means that the variety of
$\mathcal{M}_{A}$ decomposes as a  product of $\mathbb{P}^{1}$ and a smaller hidden subset model.

If $n=2$ then, up to the symmetry of the $2$-cube, there are precisely 
three distinct hidden subset models which are non-degenerate: $\{\emptyset,1,2,12\}$, $\{\emptyset,1,2\}$, $\{\emptyset, 12\}$. Their models $\mathcal{M}_A$ all parametrize the full tetrahedron
$\Delta_3$ of distributions on $2^{[2]}$ or, in algebraic terms, the whole projective space $\mathbb{P}^{3}$.

If $n=3$ then, up to symmetry of the
$3$-cube, there are precisely $19$
collections $A$ of subsets of $\{1,2,3\}$
with $2 \leq |A| \leq 7$. 
Thirteen of these $19$ models $\mathcal{M}_A$
have codimension $0$, that is, they are full-dimensional
in the simplex $\Delta_7$ of  probability distributions
on $2^{[3]}$.
One of these sets is $A = \{ \emptyset, 123 \}$
which represents 
$\mathcal{M}_A = {\rm Sec}((\mathbb{P}^1)^3)$
and hence fills $\mathbb{P}^7$.
The remaining six of the $19$ models $\mathcal{M}_A$
represent three distinct varieties.
The first of them is the hyperdeterminantal hypersurface
${\rm Tan}( (\mathbb{P}^1)^3)$. The other two varieties are
a line and a point in cumulant space:
$$\begin{array}{c|c|c}
\mbox{hidden subset model} & \mbox{variety} & \mbox{codimension}\\
\hline
\{\emptyset, 1,2,3\} \,\, {\rm or} \,\, \{\emptyset, 12, 13\}
 &  k_{123}^2 + 4 k_{12} k_{13} k_{23}  = 0 & 1\\
\{1,2\}  \, {\rm or} \,
\{\emptyset, 1,2\}  \, {\rm or} \,
\{\emptyset,1,2,12\}  
& (k_{12},k_{13},k_{23},k_{123}) = (0,0,0,*) & 3\\
\{\emptyset , 1\} & (k_{12},k_{13},k_{23},k_{123}) = (0,0,0,0) & 4\\
\end{array}
$$
The first row corresponds to non-degenerate models $\mathcal{M}_A$
that do not fill $\Delta_7$.
The situation becomes more interesting for $n \geq 4$
when we get a vast range of new models. Some of these
will be discussed and catalogued in Section~6.

\section{Context-specific independence}

This section concerns a  class of statistical models that has proved to be useful in
machine learning and computational biology \cite{Georgi15072006}, namely, 
the context-specific independence (CSI) for binary random variables.
It has been observed in \cite[\S 6.3]{oeding1} that both the tangential variety
and the secant variety of the Segre variety
are CSI models. Examples \ref{ex:tangentialrevisited} 
and \ref{ex:secantrevisited} expressed these as hidden subset models.
We here generalize this relationship by identifying
the  class of binary CSI models with a natural class of hidden subset models.

The formal specification of a CSI model is as follows. Fix a multiset of $n$ partitions $\{\pi_{1}, \pi_2, \ldots, \pi_{n}\}$
of the set $[m]=\{1,\ldots,m\}$. The model is 
\begin{equation}
\label{eq:csipara}
M(x) \,\,\, = \,\,\, \sum_{j=1}^{m} t_{j}
(1+ a_{\pi_{1}(j)} x_1)
(1+ b_{\pi_{2}(j)} x_2)  \cdots
(1+ c_{\pi_{n}(j)} x_n).
\end{equation}
Here $\pi_{j}(k)$ is the block of the $j$-th partition $\pi_j$ that contains the class $k$,
 and $t_{1},\ldots,t_{m}$ are mixing parameters for the classes. These satisfy $ t_1 + \cdots + t_m = 1$. 

If each $\pi_j$ is the partition into singletons, then we can write
$\pi_i(j) = j$ in (\ref{eq:csipara}) and the CSI model is the 
$m$-th secant variety of $(\mathbb{P}^1)^n$ in $\mathbb{P}^{2^n-1}$.
This is known in statistics as
mixture of $n$ independent binary vectors or as the {\em naive Bayes model}.
Hence every CSI model with $m$ hidden classes is a submodel of the
naive Bayes model.

The CSI model has $\sum_{i=1}^{n}|\pi_{i}|+m-1$  parameters and the dimension 
of the ambient space is $2^{n}-1$.  It is usually not identifiable, meaning
its dimension  is smaller than the number of parameters.
However, identifiability does hold for $m=2$. Here the
CSI model is the product of the Segre variety $\,(\mathbb{P}^1)^{n-k}$ and the first secant variety of $(\mathbb{P}^{1})^{k}$, where $k=\#\{j\in [n]:\pi_{j}=1|2\}$. 
This is the graphical model represented 
by a directed star tree with a hidden binary variable and $k$ leaves together with $n-k$ isolated nodes. 

From statistical point of view it is sensible to assume the following:
\begin{itemize}
\item[(A1)] All partitions in the model specification have at least two blocks. 
\item[(A2)] There is no pair of elements $\{i,j\}$ such that for every partition in the model specification both $i$ and $j$ are in the same block of this partition. 
\end{itemize}
If (A1) is violated then one random variable is independent of all others. 
Taking an appropriate margin, we can constrain our analysis to the remaining variables. If (A2) does not hold then the classes $i$ and $j$ can be joined to form a single class without changing the model. If $m=n=3$ then up to symmetry we have three CSI models satisfying (A1) and (A2). The first case is $\pi_1 = 1|23$
$\pi_2 =  2|13$, and $\pi_3 = 3|12$. This is precisely our  hyperdeterminantal hypersurface
${\rm Tan}((\mathbb{P}^1)^3) = V( k_{123}^2 + 4 k_{12} k_{13} k_{23})$.
The other two CSI models represent all distributions on $2^{[3]}$:
$$ (\pi_{1}=1|23,\, \pi_{2}=12|3,\,\pi_{3}=1|2|3)\quad\mbox{or} \quad
(\pi_{1}=1|23,\, \pi_{2}=\pi_{3}=1|2|3). $$

In the remainder of this section we study a special class of CSI models.
Namely, we shall require that
each partition $\pi_i$ has precisely two blocks. 
  We call these models the {\em CSI split models}. 
  Thus a CSI split model is represented by a collection
  $\{\pi_1,\pi_2, \ldots ,\pi_n\}$ of splits  of the
  set $[m]$ of hidden states.
  The following result identifies these models with the models  in Section~5.
    
\begin{prop} \label{prop:csi=hsm}
The CSI split models are precisely the hidden subset models.
\end{prop}

\begin{proof}
Let $\mathcal{M}_A$ be the hidden subset model defined by $A=\{J_{1},\ldots, J_{m}\} \subseteq 2^{[n]}$.
This is written as a CSI split model with $m$ hidden classes
by taking the $n$ partitions $\pi_{1},\ldots, \pi_{n}$ of $[m]$ to be $\pi_{i}=I|I^{c}$, 
where $\ell \in I$ whenever  $i \in J_\ell $. Conversely, suppose we are
given a CSI split model  $\{\pi_{1},\ldots, \pi_{n}\}$. Then we regard
$\pi_i$ as an ordered partition, and we recover the $m$ subsets in $A$ by 
taking $J_\ell$ to be the set of all $i \in [n]$ such that
$\ell$ is in the first part of $\pi_i$.
These transformations lead to identical parametrizations,
and hence the corresponding models in $\Delta_{2^n-1}$ coincide.
\end{proof}

We classified all hidden subset models and hence all CSI split models for $n=3$
in the end of the previous section. The next case $n=4$ is much more interesting, as it offers
a considerably wider range of possibilities. The classification for $n=4$ will occupy us in the
rest of this section. 
%We begin with an illuminating example.

\begin{exmp}[The hyperdeterminant as CSI model] \label{ex:hypercsi}
Let $n= 4$,  $m=7$, and consider the hidden subset model $\mathcal{M}_A$ where
$A = \{ \emptyset, 12, 13, 14, 23, 24, 34 \}$. The corresponding CSI split
model is
$\,\{ 1234|567, 1256|347, 1357|246, 1467|235\}$. In algebraic geometry, the model $M_A$ corresponds to the
{\em second osculating variety} of the Segre variety $(\mathbb{P}^1)^4$.
This is a hypersurface of degree $24$ in $\mathbb{P}^{15}$, namely, it is the hypersurface defined by
 the {\em $2 {\times} 2 {\times} 2 {\times} 2$-hyperdeterminant}.
This result was pointed out to us by Luke Oeding and Giorgio Ottaviani, and we
can easily verify it by a direct computation. The fact that ${\rm codim}(\mathcal{M}_A) = 1$
is verified by computing the rank of the Jacobian of
the parametrization (\ref{eq:remarkable}) for random parameter values.
The fact that $\mathcal{M}_A$ equals $\{{\rm Det}(P) = 0\}$
is verified by plugging the parametrization (\ref{eq:remarkable})  into the formula
with $13819$ monomials found in  Theorem \ref{thm:hyperdet}. 
We note that this model remains the same if we augment
 $A$ to $\{\emptyset,1,2,3,4, 12, 13, 14, 23, 24, 34\}$.
\qed
\end{exmp}

We now come to the classification of CSI split models for $n= 4$.
Each model lives in the space $\mathbb{C}^{11}$ with coordinates 
$k_{12},\ldots,k_{34},k_{123},\ldots,k_{234},k_{1234}$.

\begin{prop}\label{thm:listcsi4}
Up to symmetry, for $n{=}4$, there are
% precisely
 $380$ CSI split models satisfying (A1) and (A2). The number of models with $m$ hidden classes is
$$
0 , 1 , 3 , 13 , 24 , 47 , 55 , 73 , 56 , 50 , 27 , 19 , 6 , 4 , 1 , 1 
\qquad \hbox{for} \quad m = 1,2,\ldots,16.
$$
In Table \ref{tab:csi4} we list 
the codimension and degree for all $\,17$ models for $m\leq 4$.
\end{prop}

%\begin{table}[b]
\begin{center}
\begin{table}
%\begin{small}
$$\begin{array}{c|c|c|c}
\mbox{Hidden subset model}  & \mbox{CSI split model}& \mbox{codimension} & \mbox{degree}\\
%\hline
%\{\emptyset, 1,2,12,34\} & 135|24,125|34,1234|5,1234|5 & 5 &  44\\
%\{\emptyset, 1,13,14,234\} & 15|234,1234|5,124|35,1234|45 & 3 &  \\
%\{\emptyset, 1,12,34,123\} & 14|235,124|35,123|45,1235|4 & 4 &  \\
%\{\emptyset, 1,2,34,123\} & 134|25,124|35,123|45,1235|4 & 4 &  \\
%\{\emptyset, 1,12,13,34\} & 15|234,1245|3,123|45,1234|5 & 4 &  \\
%\{\emptyset, 1,12,23,34\} & 145|23,125|34,123|45,1234|5 & 4 &  \\
%\{\emptyset, 2,3,12,14\} & 123|45,135|24,1245|3,1234|5 & 5 & \\
%\{\emptyset, 3,12,14,123\} & 12|345,124|35,134|25,1235|4 & 5& \\
%\{\emptyset, 12,13,14,34\} & 15|234,1345|2,124|35,123|45 & 3& \\
%\{\emptyset, 12,13,14,123\} & 1|2345,134|25,124|35,1235|4 & 6& \\
%\{\emptyset, 3,12,14,124\} & 12|345,124|35,1345|2,123|45 & 5& \\
%\{\emptyset, 1,2,13,24\} & 135|24,124|35,1235|4,1234|5 & 5& \\
%\{\emptyset, 1,2,13,14\} & 13|245,1245|3,1235|4,1234|5 & 6 & 29 \\
%\{\emptyset, 2,3,14,123\} & 123|45,134|25,124|35,1235|4 & 4 & \\
%\{\emptyset, 1,12,34,134\} & 145|23,1245|3,123|45,123|45 & 5 & \\
%\{\emptyset, 1,12,13,14\} & 1|2345,1245|3,1235|4,1234|5 & 7 & \\
%\{\emptyset, 12,13,14,234\} & 15|234,134|25,124|35,123|45 & 3& \\
%\{\emptyset, 12,34,134,234\} & 135|24,134|25,124|35,12|345 & 4& \\
%\{\emptyset, 2,13,14,123\} & 12|345,134|25,124|35,1235|4 & 3& \\
%\{\emptyset, 12,23,34,134\} & 134|25,145|35,12|345,123|45 & 3& \\
%\{\emptyset, 1,14,23,234\} & 145|23,123|45,123|45,124|35 & 5& \\
%\{\emptyset, 1,2,134,234\} & 135|24,124|35,123|45,123|45 & 4& \\
%\{\emptyset, 1,23,124,234\} & 135|24,12|345,124|35,123|45 & 4& \\
%\{\emptyset, 12,13,24,34\} & 145|23,135|24,124|35,123|45 & 3& \\
\hline
\{\emptyset, 12,13,14\} & 1|234,2|134, 3|124,4|123& 7  & 20 \\
\{\emptyset, 12,13,4\} & 14|23,2|134, 3|124,4|123& 6 & 29\\
\{\emptyset,1,2,34\} & 2|134,3|124,4|123,4|123& 6 & 29\\
\{\emptyset,1,23,234\} & 2|134,34|12,34|12,4|123& 6 & 23\\
\{\emptyset,1,234,1234\} & 24|13,34|12,34|12,34|123& 6 & 23\\
\{\emptyset, 1,2,134\} & 13|24,3|124, 4|123,4|123 & 5& 44\\
\{\emptyset,1,12,234\} & 23|14,34|12,4|123,4|123& 5 & 44\\
\{\emptyset,1,123,234\} & 23|14,34|12,34|12,4|123& 5 & 44\\
\{\emptyset,1,23,124\} & 24|13,34|12,3|124,4|123& 5 & 31\\
\{\emptyset,12,134,234\} & 23|14,24|13,34|12,34|12& 5 & 22\\\{\emptyset, 12,13,24\} & 14|23,13|24, 3|124,4|123& 4 & 44 \\
\{\emptyset, 13,23,124\} & 13|24,12|34, 14|23,4|123&  4 & 38\\
\{\emptyset,12,34,1234\} & 24|13,24|13,34|12,34|12& 4 & 11\\
\hline
\{\emptyset,1,234\} & 2|13,3|12,3|12,3|12& 6& 23 \\
\{\emptyset, 12,134\} & 1|23,2|13, 3|12,3|12& 6& 29 \\
\{\emptyset,12,34\} & 2|13,2|13,3|12,3|12& 5& 44 \\
\hline
\{\emptyset,1234\} & 1|2,1|2,1|2,1|2 & 6  & 23 \\
\end{array}
$$
\caption{
The~$17$~non-degenerate CSI split models
on $n{=}4$ binary variables with $m \leq 4 $ hidden classes, up to symmetry.}
\label{tab:csi4}
%\end{small}
\end{table}
\end{center}

For our classification we used the representation 
of each CSI split model as a hidden subset model $\mathcal{M}_A$, 
where $A \subseteq 2^{[4]}$, given by Proposition \ref{prop:csi=hsm}.
A choice of $A$ is also displayed for each model in
Table \ref{tab:csi4}.  Note that two distinct models may define the same variety.
By symmetry we can assume $\emptyset \in A$.
We first generated the list  of all non-degenerate sets $A$ of subsets of $\{1,2,3,4\}$ 
containing $\emptyset$, we then computed orbits under the symmetry group
of the $4$-cube, and finally we selected one representative per orbit.
To compute the codimension $c$ of $\mathcal{M}_A$,
we evaluated the rank of the Jacobian of the polynomial map (\ref{eq:remarkable})
at random values of the parameters.
By {\em degree} in Proposition \ref{thm:listcsi4} we mean  the number
of complex solutions on $\mathcal{M}_A$ of a system
of $11-c$ inhomogeneous linear equations
with random coefficients in the $11$ unknowns $k_I$.
We used  {\tt Macaulay2} \cite{M2} to count the
number of solutions to these equations.
It would be desirable to compute the defining prime
ideals for all models in Proposition~\ref{thm:listcsi4},
but we found this to be difficult for $m \geq 5$.
The $380$ models represent a nice suite of test problems
for {\em implicitization} in computer algebra.
We close the section with one easy instance.

\begin{exmp}
The model $\mathcal{M}_A$ with $A = \{\emptyset , 12, 34, 1234 \}$ has
CSI representation $\{12|34, 12|34, 13|24, 13|24\}$.
Its prime ideal in cumulant coordinates is
$$ 
\begin{small}
\begin{matrix}
\langle k_{13} k_{24} - k_{14} k_{23},\,\,
k_{13} k_{124} - k_{14} k_{123},\,\,
k_{13} k_{234} - k_{23} k_{134},  \,\,
k_{14} k_{234} - k_{24} k_{134}, \\
k_{23} k_{124} - k_{24} k_{123}, \,
k_{23} k_{1234} - k_{234} k_{123} + 2 k_{14} k_{23}^2, \,
k_{13} k_{1234} - k_{134} k_{123} + 2 k_{14} k_{13} k_{23}, \\
k_{23} k_{1234} - k_{234} k_{124} + 2 k_{14} k_{24} k_{23}, \,\,\hbox{and} \,\,\,
k_{14} k_{1234} - k_{134} k_{124} + 2 k_{14}^2 k_{23} \rangle.
\end{matrix}
\end{small}
$$
This CSI split model has codimension $4$ and degree $11$.
 \qed
\end{exmp}

\section{Semialgebraic geometry of the space of cumulants}\label{sec:semialg}

In the previous sections we studied binary cumulant varieties as objects of complex algebraic geometry.
We examined their dimension, parameterization, and defining prime ideal, but we largely ignored
the issue that parameters and probabilities are real and non-negative.
In statistical applications, however, it is essential
to work over the real numbers and to pay attention to the pertinent inequalities.
In this section seek to address this omission by asking
 the following fundamental question:
Which $2 {\times} 2 {\times} \cdots {\times} 2$-tables $K = (k_I)_{I \subseteq [n]} $ 
with entries in the real numbers represent
the cumulants of actual probability distributions
$P = (p_I)_{I \subseteq [n]}$?

Our object of study is
the image of the polynomial map
$\Delta_{2^n-1} \rightarrow \mathbb{R}^{2^n-1}$
taking probability distributions $P$ to their 
cumulants $K$. This image is 
denoted $\mathcal{K}_n$. We call it
the {\em space of cumulants}.
 The space of cumulants $\mathcal{K}_n$ is a 
semi-algebraic subset of $\mathbb{R}^{2^n-1}$. This
means that it has a description in terms of polynomial
inequalities in the $k_I$. We begin by offering
a convenient representation of these inequalities.

\begin{prop} \label{prop:rhojk}
The space of cumulants $\mathcal{K}_n$
is a basic semialgebraic set in $\mathbb{R}^{2^n-1}$.
It consists of the solutions of the polynomial inequalities
\begin{equation}
\label{eq:rhojk}
  \sum_{\pi\in \Pi([n])}\prod_{B\in \pi}
\rho_J(k_{B}) \,\,\geq \,\, 0\qquad\mbox{for all } \, J\subseteq [n].
\end{equation}
\end{prop}

\begin{proof}
The set $\mathcal{K}_n$ being {\em basic semialgebraic}
means that it is described by a finite conjunction of polynomial 
inequalities. That conjunction is (\ref{eq:rhojk}),
and we shall now prove it.
The moment $\mu_{1\cdots n}$ agrees with the probability $p_{1\cdots n}$,
so it is non-negative on $\mathcal{K}_n$. Expressing $\mu_{1\cdots n}$
in terms of cumulants as in (\ref{eq:momincum}), 
$$  p_{1\cdots n} \,\,\,  = \,  \sum_{\pi\in \Pi([n])}\prod_{B\in \pi} k_{B} \,\,\,\, \geq \,\,\,\, 0. $$
By applying the transformation $\rho_J$ from (\ref{actionensures}) to this inequality, 
we obtain $p_{J^c}  = \rho_{J}(p_{1\cdots n}) \geq 0$. This translates into
the inequality (\ref{eq:rhojk}) in cumulants. Since the transformation
$P \mapsto K$ is invertible, we see that $\mathcal{K}_n$
has the desired~representation.
\end{proof}

\begin{figure}
\includegraphics[scale=0.6]{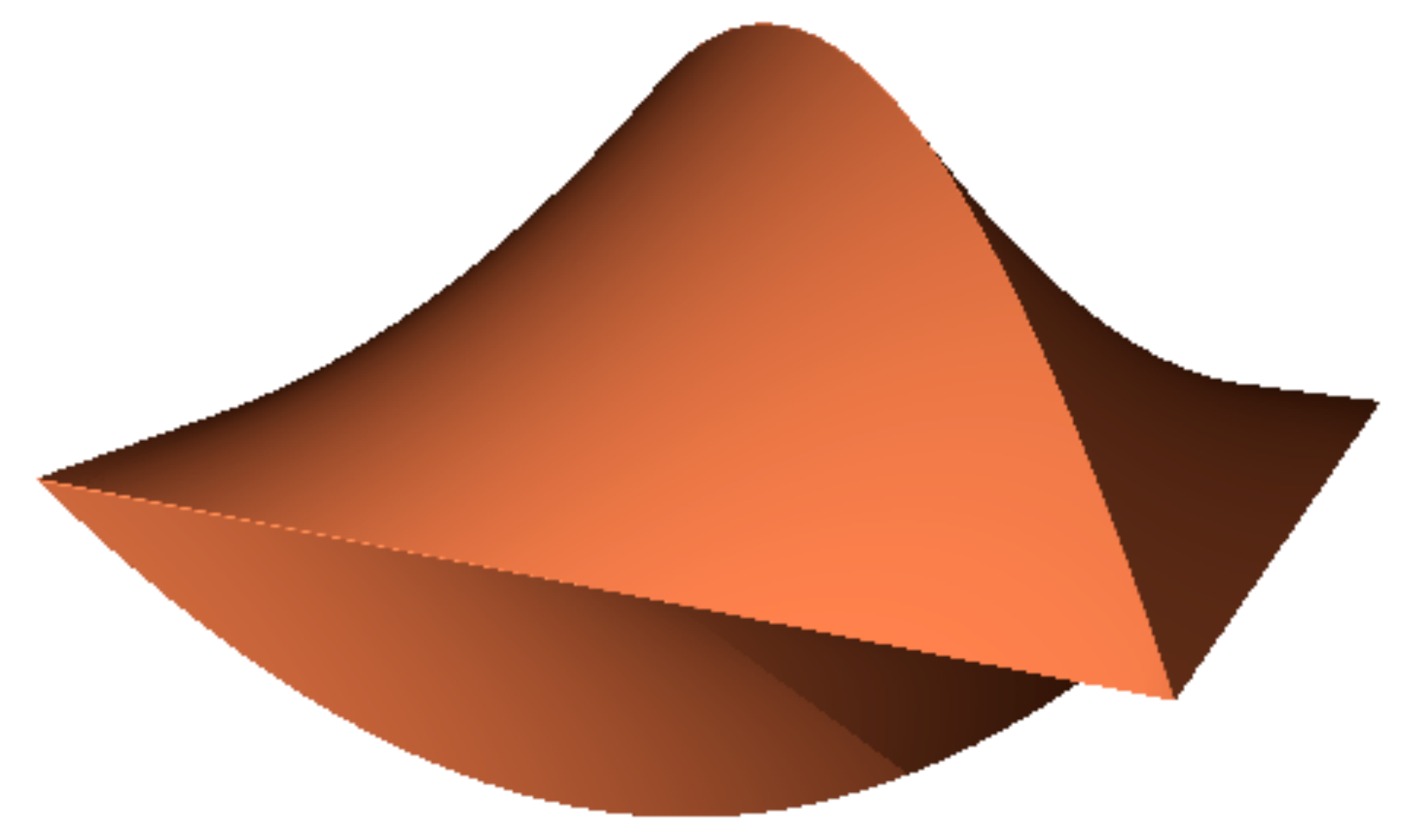}
\caption{The space of cumulants $\mathcal{K}_2$
is the solution set of (\ref{eq:constraints})}\label{fig:pillow}
\end{figure}

\begin{exmp}[Space of cumulants for $n=2$] \label{eq:Nzwei}
The probability distributions $P$ on the subsets of $\{1,2\}$ form 
a tetrahedron, and we map this tetrahedron into
the $3$-space with coordinates $(k_1,k_2,k_{12})$.
The image of this map is the space of cumulants $\mathcal{K}_2$.
Proposition \ref{prop:rhojk} gives the semialgebraic representation:
\begin{equation}
\label{eq:constraints}
\begin{matrix}
\hbox{Inequalities defining $\Delta_3$} &  & \hbox{Inequalities defining $\mathcal{K}_2$}  \\
p_{12} \geq 0 & & k_{12} \,\geq\, -k_1 k_2,\\
 p_1 \geq 0 &  & k_{12} \,\leq \,k_2(1-k_1),\\
 p_2 \geq 0 &  & k_{12}\,\leq\, k_1(1-k_2) ,\\
 p_\emptyset \geq  0 & & \,\,\,\, k_{12}\,\geq\, -(1-k_1)(1-k_2).
 \end{matrix}
\end{equation}
The solution set of these four quadratic
inequalities is depicted in Figure \ref{fig:pillow}.
In this diagram we see clearly how $\mathcal{K}_2$ arises
as a non-linear image of the tetrahedron $\Delta_3$.
Note that the body $\mathcal{K}_2$ is not convex.
The square $\{0 \leq k_1, k_2 \leq 1\}$ in the 
plane $\{k_{12} = 0\}$ is the image of the independence
model $\{p_\emptyset p_{12} = p_1 p_2 \}$, which contains four of the six edges of $\Delta_3$.
The other two edges of the tetrahedron $\Delta_3$
are the quadratic curves that form the ridges at the
top and the bottom of the $\mathcal{K}_2$.
\qed 
\end{exmp}

\begin{exmp}[The space of cumulants for $n=3$]
We now consider the simplex $\Delta_7$
of distributions on subsets of $\{1,2,3\}$. Its image
in cumulant coordinates is the 
$7$-dimensional closed basic semialgebraic set
$\mathcal{K}_3$. Both
$\Delta_7$ and $\mathcal{K}_3$ are defined 
by the constraints that the following eight
expressions should be non-negative:
\begin{small}
$$
\begin{matrix}
p_{123} \!\! & \!\! = \!\! & \!\!
 \mu_{123} \!\! & \!\! = \!\! & \!\!  k_{123} + k_{12} k_3 + k_{13} k_2 + k_{23} k_1 + k_1 k_2 k_3 \\
p_{12} \!\! & \!\! = \!\! & \!\!
  -\mu_{123} + \mu_{12} \!\! & \!\! = \!\! & \!\! -k_{123}-k_{12} (k_3-1) -k_{13} k_2
-k_{23} k_1-k_1 (k_2-1) k_3 \\
p_{13} \!\! & \!\! = \!\! & \!\!
  -\mu_{123} + \mu_{13} \!\! & \!\! = \!\! & \!\! -k_{123} - k_{12} k_3
-k_{13} (k_2-1) - k_{23} k_1 - (k_1-1) k_2 k_3 \\
p_{23} \!\! & \!\! = \!\! & \!\!
  -\mu_{123} + \mu_{23} \!\! & \!\! = \!\! & \!\! -k_{123} - k_{12} k_3 - k_{13} k_2
-k_{23} (k_1-1)-(k_1-1) k_2 k_3 \\
\end{matrix}
$$
$$
\begin{matrix}
p_{1} \!\! & \!\! = \!\! & \!\!
   \mu_{123} {-} \mu_{12} {-}\mu_{13} {+} \mu_1 \!\! & \!\! = \!\! & \!\! 
k_{123} {+} k_{12} (k_3{-}1) {+} k_{13} (k_2{-}1) {+} k_{23} k_1 {+} k_1 (k_2{-}1) (k_3{-}1) \\
p_{2} \!\! & \!\! = \!\! & \!\!
  \mu_{123} {-} \mu_{12} {-} \mu_{23} {+} \mu_2 \!\! & \!\! = \!\! & \!\! k_{123} {+} k_{12} 
(k_3{-}1) {+} k_{13} k_2 {+} k_{23} (k_1{-}1) {+} (k_1{-}1) k_2 (k_3{-}1) \\
p_{3} \!\! & \!\! = \!\! & \!\!
   \mu_{123}{-}\mu_{13}{-}\mu_{23} {+} \mu_3 \!\! & \!\! = \!\! & \!\! k_{123} {+} k_{12} k_3 
{+} k_{13} (k_2{-}1) {+} k_{23} (k_1{-}1) {+} (k_1{-}1) (k_2{-}1) k_3 
\end{matrix}
$$
\vskip -0.34cm
$$
\begin{matrix}
\quad p_\emptyset & = &
 -\mu_{123} + \mu_{12} + \mu_{13} + \mu_{23} - \mu_1-\mu_2-\mu_3 + 1 \\ & = & 
 -k_{123} - k_{12} (k_3 {-} 1) - k_{13} (k_2 {-} 1) - k_{23} (k_1 {-} 1) - (k_1 {-} 1) (k_2 {-} 1) (k_3 {-} 1)
\end{matrix} \qquad
$$
\end{small}
Thus the space $\mathcal{K}_3$ is defined by eight cubic inequalities in $\mathbb{R}^7$. \qed
\end{exmp}

Equipped with the inequality description of $\mathcal{K}_n$ we can now
try to answer questions about the geometry of cumulants of probability
distributions. One natural such question is to identify the smallest
box containing $\mathcal{K}_n$. This is equivalent to find a tight
upper and lower bound on  the possible values of the cumulants $k_I$.
The following problem was suggested by
 Gian-Carlo Rota and his collaborators in \cite{bruno1999probability}:
 \begin{equation}
\label{eq:rotaopt}
{\rm Maximize} \,\, | k_{12\cdots n}| \,\,\, \hbox{subject to} \,\,\, k \in \mathcal{K}_n.
\end{equation}

In this problem, the absolute value sign around $k_{12 \cdots n}$
can be removed because
$k \in \mathcal{K}_n$ implies $-k \in \mathcal{K}_n$. 
This is shown in \cite[Proposition 3.1]{bruno1999probability},
and it also follows directly from the symmetries 
 in Corollary \ref{cor:sym-dif}.
Let ${\bf k}^*_n$ denote the optimal value of (\ref{eq:rotaopt}).
Figure \ref{fig:pillow} shows that ${\bf k}^*_2 = 1/4$,
and one easily derives an algebaic proof from the inequalities
in Example \ref{eq:Nzwei}. The probability distribution
$\,p_\emptyset = p_{12}  = \frac{1}{2}$ attains $\,{\bf k}^*_2 = 1/4$.
 It has been conjectured by Bruno, Rota and Torney
 \cite{bruno1999probability} that the analogous 
 distribution solves the optimization problem (\ref{eq:rotaopt})
 for all even values of~$n$:

\begin{conj}\label{conj:rota} \cite[bottom of page 16]{bruno1999probability} \ \
If $n \geq 2$ is an even integer then
$$ {\bf k}^*_n \,\,= \,\,
\kappa_n (\frac{1}{2}) \,\, = \,\,(-1)^{n/2} \sum_{i=1}^n \, (-\frac{1}{2})^i  \cdot \gamma_{i, n} . $$
This value is attained by the
probability distribution  $\,p_{\emptyset}=p_{1\cdots n}=\frac{1}{2}$.
\end{conj}

Here $\kappa_n$ is the polynomial in  (\ref{eq:necklaces}).
The first values of the bound $\kappa_n(\frac{1}{2})$ are
$\frac{1}{4}, \frac{1}{8}, \frac{1}{4}, \frac{17}{16}, \frac{31}{4}$
for $n = 2,4,6,8,10$.  It has been remarked in \cite{bruno1999probability} that
$$
\qquad
\bigl|\kappa_n(\frac{1}{2}) \bigr| \,\,\sim\,\, 2\frac{1}{\pi^{2n}}(2n-1) ! \qquad
\hbox{for} \,\,\, n \gg 0 . 
$$
If $n \geq 3$ is an odd integer then 
$\kappa_n(\frac{1}{2}) = 0$,
and no conjectured value
for ${\bf k}^*_n$ has been suggested in \cite{bruno1999probability}.
Using recent computational advances in certified polynomial optimization, 
we attacked the problem (\ref{eq:rotaopt}) for $n=3$ and $n=4$,
thus confirming the conjecture of Bruno, Rota and Torney   in the first non-trivial case.
Namely, we found that
the upper bound on cumulants of probability set functions satisfies
\begin{equation}
\label{eq:k*34}
 {\bf k}^*_3 \,\,\,= \,\,\, {\bf k}^*_4 \,\,\, = \,\,\, \frac{1}{8}.
 \end{equation}

For $n=3$ we used
the software {\tt Bermeja} \cite{bermeja}
to compute a  {\em sums of squares} certificate
via semidefinite programming. 
We are grateful for the help provided by Philipp Rostalski.
Let us now explain this certificate.
We consider the following cubic polynomial
in the seven moment coordinates~$\mu_I$:
$$
  \frac{1}{8}- k_{123}  \quad = \quad
   \frac{1}{8} - \mu_{123} + \mu_1 \mu_{23} + \mu_2 \mu_{12} + \mu_3 \mu_{12}  - 2 \mu_1 \mu_2 \mu_3.  $$
Our aim is to prove that this polynomial is non-negative on the simplex $\Delta_7$.
We do this by rewriting the polynomial in the following special form
\begin{small}
\begin{equation}
\label{eq:mustcheck1} \,\,
  \frac{1}{8} - k_{123}  \,\, = \,\,
\sigma_\emptyset
+ \sigma_1 \mu_1 + \sigma_2 \mu_2 + \sigma_3 \mu_3
+ \sigma_{12} \mu_{12} + \sigma_{13} \mu_{13} + \sigma_{23} \mu_{23}
+ \sigma_{123} \mu_{123} ,
\end{equation}
\end{small}
where each of the eight multipliers $\sigma_I$ is a 
sum of squares of linear polynomial in the moments $\mu_J$.
Each such sum of squares corresponds to a positive semidefinite quadratic form,
and it can be represented by a symmetric $8 \times 8$-matrix $\Sigma_I$ as follows:
\begin{equation}
\label{eq:mustcheck2}
\sigma_I \,\,= \,\,\mu \cdot \Sigma_I \cdot \mu^T \qquad \hbox{where}
\quad \mu \,=\, (1,\mu_1,\mu_2,\mu_3 ,\mu_{12} ,\mu_{13} ,\mu_{23} ,\mu_{123}).
\end{equation}
Our certificate for ${\bf k}^*_3 = 1/8$ is a tuple
$\bigl(\Sigma_\emptyset,\Sigma_1,\Sigma_2,\Sigma_3 ,\Sigma_{12} 
,\Sigma_{13} ,\Sigma_{23} ,\Sigma_{123} \bigr)$
of positive semidefinite symmetric $8 {\times} 8$-matrices such that
(\ref{eq:mustcheck1}) and (\ref{eq:mustcheck2}) hold.
Finding such a tuple of matrices is an instance of semidefinite programming.

We attempted to find a similar proof
for the second identity $\mathbf{k}^*_4 = 1/8$
but the computations required turned out to be too difficult so far.
The idea was to take advantage of the symmetries
preserves the optimization problem (\ref{eq:rotaopt}). This is a group
of order $192$, and has index $2$ in the symmetry group of
the $4$-cube. Our hope was to use the
the dual moment formulation due
to Riener~{\em et al.} in \cite{JLRT}, but this did yet terminate successfully.
Instead, we verified the identify $\mathbf{k}^*_4 = 1/8$ by running 
numerous applications of standard implementations of numerical
optimization in {\tt R} and {\tt Matlab}. Running these hill climbing methods
from a multitude of different starting values verifies the desired result with
very high confidence.

\bibliographystyle{siam}
%\bibliography{../!bibliografie/algebraic_statistics}

\end{document}